\documentclass{elsarticle}

\usepackage{graphics} 
\usepackage{epsfig}   
\usepackage{mathptmx} 
\usepackage{times}    
\usepackage{amsmath}  
\usepackage{amssymb}  
\usepackage{subfig}
\usepackage{booktabs}
\usepackage{verbatim}

\newcommand{\hF}{\mathbf{F}}
\newcommand{\hG}{\mathbf{G}}
\newcommand{\hH}{\mathbf{H}}
\newcommand{\hx}{\mathbf{x}}
\newcommand{\hu}{\mathbf{u}}
\newcommand{\hw}{\mathbf{w}}

\newcommand{\dd}{\mathrm{d}}
\newcommand{\opt}{\lambda}

\newtheorem{theorem}{Theorem}

\begin{document}

\begin{frontmatter}

\title{A relaxation technique to ensure feasibility in stochastic control with input and state constraints\tnoteref{thx1}}

\tnotetext[thx1]{The authors gratefully acknowledge financial support by the European Commission project UnCoVerCPS under grant number 643921.}

\author[PoliMi]{Luca Deori}\ead{luca.deori@polimi.it}
\author[PoliMi]{Simone Garatti}\ead{simone.garatti@polimi.it}
\author[PoliMi]{Maria Prandini}\ead{maria.prandini@polimi.it}

\address[PoliMi]{Dipartimento di Elettronica, Informazione e Bioingegneria, \\ Politecnico di Milano, Via Ponzio 34/5, 20133 Milano, Italy}

\begin{abstract}
We consider a stochastic linear system and address the design of a finite horizon control policy that is optimal according to some average cost criterion and accounts also for probabilistic constraints on both the input and state variables. This finite horizon control problem formulation is quite common in the literature and has potential for being implemented in a receding horizon fashion according to the model predictive control strategy. Such a possibility, however, is hampered by the fact that, if the disturbance has unbounded support, a feasibility issue may arise. 
In this paper, we address this issue by introducing a constraint relaxation that is effective only when the original problem turns out to be unfeasible and, in that case, recovers feasibility as quickly as possible. This is obtained via a cascade of two probabilistically-constrained optimization problems, which are solved here through a computationally tractable scenario-based scheme, providing an approximate solution that satisfies the original probabilistic constraints of the cascade, with high confidence. A simulation example showing the effectiveness of the proposed approach concludes the paper.
\end{abstract}


\end{frontmatter}

\section{Introduction}

This paper deals with the problem of designing a finite-horizon optimal control policy for a discrete-time stochastic linear system subject to constraints on both the state and the input. Specifically, we consider the system
\begin{align}\label{eq:sys}
x_{t+1}=Ax_t+Bu_t+B_ww_t,
\end{align}
where  $x_t \in \mathbb{R}^n$ is the state, $u_t \in \mathbb{R}^m$ is the control input and $w_t \in \mathbb{R}^{n_w}$ is an additive stochastic disturbance. Matrices $A$, $B$, and $B_w$ have appropriate dimensions so as to make \eqref{eq:sys} consistent. The probability distribution of $w_t$, say $\mathbb{P}$, is assumed to be known and it may have an unbounded support.  Without loss of generality, we assume that $n_w\le n$ and $B_w$ is full column rank. The state is accessible, i.e., at every $t$ a noise-free measurement of $x_t$ becomes available.

The following disturbance feedback parametrization for the control input is adopted:
\begin{align}\label{eq:input_param}
&u_{t}=\gamma_t + \sum_{\tau=0}^{t-1}\theta_{t,\tau}w_{\tau},
\end{align}
where $\gamma_t \in \mathbb{R}^{m}$ represent open-loop terms, while $\theta_{t,\tau} \in \mathbb{R}^{m\times n_w}$ are the disturbance feedback gains. Both $\gamma_t$ and $\theta_{t,\tau}$ are degrees of freedom to be chosen.  Note that the stochastic disturbance $w_\tau$ appearing in it can be recovered from the measurements of the state according to
\begin{align}\label{eq:reco_W}
 w_\tau=B^{\dag}_{w}(x_{\tau+1} - A x_\tau - B u_{\tau}),
\end{align}
where $B^{\dag}_w$ denotes the pseudo-inverse of $B_w$. This expression reveals that the disturbance feedback control policy in \eqref{eq:input_param} is in fact a state feedback control policy. Parametrization \eqref{eq:input_param} was first proposed in \cite{GKM2006}, where it was shown that the family of policies in \eqref{eq:input_param} is indeed equivalent to the family of affine state feedback policies $u_t = \bar{u}_t + K_t [x_t^T\ x_{t-1}^T \ \dots \ x_0^T]^T$. To be precise, for every choice of $\bar{u}_t , K_t$ there exists a parametrization $\gamma_t , \theta_{t,\tau}$ in \eqref{eq:input_param} returning the same control action, and viceversa. The great advantage of \eqref{eq:input_param} is that, differently from other parametrizations, the input $u_t$ and the state $x_t$ are affine functions of the design parameters $\gamma_t$ and $\theta_{t,\tau}$, which yields clear computational benefits.

The objective is to design the parameters $\gamma_t$ and $\theta_{t,\tau}$ so as to minimize a cost function over a finite  time horizon of length $M$, while accounting for constraints on the input and state variables. This problem may arise per-se in some applications (for instance, the positioning of the end-effector of an industrial robot equipped with a robot re-initialization device), but its significance mainly lies in the fact that it can be adopted in a Model Predictive Control (MPC) scheme, where it is repeatedly solved at every time step, \cite{maciejowski_2001,Camacho_Bordons_2004,Rawlings_Mayne_2009,kouvaritakis2016model}.

In our formulation, we admit as cost any strictly convex function $J$ of the parameters $\gamma_t$ and $\theta_{t,\tau}$ over the horizon $0,1,\ldots,M-1$. Plainly, the most common situations is when $J$ is defined as a function of the input and the state. A typical choice is the average quadratic cost
\begin{align}\label{eq:cost}
J&=\mathbb{E}\left[\sum_{t=1}^{M} x_t^TQx_t + \sum_{t=0}^{M-1} u_t^TRu_t \right],
\end{align}
where $Q$ and $R$ are symmetric and positive semi-definite matrices, and $\mathbb{E}$ denotes expectation with respect to the probability distribution $\mathbb{P}$. In this case, a sufficient condition for strict convexity to hold is that matrices $R$ and $\mathbb{E}[\hw \hw^T]$ are positive definite, see \cite{Franklin16}. Other choices are however possible.

As for the input and state constraints, we assume that they are expressed as
\begin{equation} \label{eq:ux_constr}
f(u_0, \dots,u_{M-1})  \le 0 \ \wedge \ g(x_1, \dots,x_{M},u_0, \dots,u_{M-1}) \le 0
\end{equation}
where $\wedge$ stands for ``and'', $f: \mathbb{R}^{mM} \rightarrow \mathbb{R}^{p_u}$ and $g: \mathbb{R}^{(n+m)M} \rightarrow \mathbb{R}^{p_y}$ are continuous convex vector-valued functions, and the inequalities are meant component-wise.  For example, a typical requirement is that the norm of the input and of some output variable are kept within an admissible range, in which case we have
\begin{align*}
& f(u_0, \dots,u_{M-1})=\begin{bmatrix}
\|u_{0}\| \\
\vdots\\
\|u_{M-1}\| \\
\end{bmatrix}
-\bar u  & g(x_1, \dots,x_{M},u_0, \dots,u_{M-1})= \begin{bmatrix}
\|Cx_{1}\| \\
\vdots\\
\|Cx_{M}\| \\
\end{bmatrix}
-\bar y,
\end{align*}
where the vectors $\bar u$ and $\bar y$ defines the maximum allowed magnitude at each time instant and $\| \cdot \|$ denotes some norm of interest. Note that $g$ allows for joint state and input constraints along the temporal horizon of interest and the constraints expressed by $f$ could be incorporated in $g$. To ease further explanations, we however keep the constraints that depend on the input only separate from the others.

It should be noted that constraints \eqref{eq:ux_constr} cannot be directly imposed since they miss to specify how to account for the presence of the stochastic disturbance affecting both the state and the input variables. Since the disturbance support is possibly unbounded, we assume that constraints are enforced probabilistically, namely, constraints \eqref{eq:ux_constr} are required to hold with a certain (usually high) probability $1-\epsilon$, where $\epsilon\in (0,1)$ is a user-chosen parameter:
\begin{align}\label{eq:Pconstr_y}
\mathbb{P}\{&f(u_0, \dots,u_{M-1})  \le 0\ \wedge g(x_1, \dots,x_{M},u_0, \dots,u_{M-1}) \le 0  \} \ge 1-\epsilon.
\end{align}
This probabilistic formulation of constraints has actually become common in the recent literature on constrained stochastic control, \cite{BB2007,Ono2008,Primbs_2007,Primbs_2009,BO2009,CACL2011,cannon2011stochastic}. Altogether, the optimal design problem we are considering is as follows:
\begin{align}\label{eq:prob_original}
&\min_{\gamma_i, \theta_{i,j}}\ J \\
&\quad \text{subject to: }  \mathbb{P}\{f(u_0, \dots,u_{M-1})  \le 0\ \land g(x_1, \dots,x_{M},u_0, \dots,u_{M-1}) \le 0  \} \ge 1-\epsilon. \nonumber
\end{align}

Although \eqref{eq:prob_original} is the most natural formulation for many problems of interest, a feasibility issue arises, precisely because of the presence of the requirement on $g(x_1, \dots,x_{M}, \allowbreak u_0, \dots,u_{M-1})$. As a matter of fact, the stochastic disturbance $w_t$ enters additively the system dynamics, and, since the input $u_t$ depends on the disturbance up to time $t-1$ at most, the dependence of $x_{t+1}$ on $w_t$ cannot be canceled. Since $w_t$ has possibly unbounded support and given the limitation imposed by the system dynamics and by the constraints on the input variable, it may then be that, depending on the system initialization $x_0$, no choice of $\gamma_t,\theta_{t,\tau}$ exists such that $g(x_1, \dots,x_{M},u_0, \dots,u_{M-1}) \le 0$ is attained with the required probabilistic level.
{It is perhaps worth noticing that probabilistic constraints of the kind $\mathbb{P} \{ f(u_0, \dots,u_{M-1})  \le 0 \} \ge 1-\epsilon$, where $g$ is not present, are instead always feasible, because, if needed, the disturbance feedback gains can be set to zero, which makes $u_t$ deterministic.}

The feasibility problem here discussed is severe because in many cases the designer has no direct control on the system initialization, which is indeed determined by exogenous causes. For example, in an MPC scheme where the optimization problem \eqref{eq:prob_original} is continuously repeated at each time step over a receding horizon and only the first calculated control action is actually implemented, the system initialization for a given problem is determined by the solutions at previous steps. At these previous steps, however, since constraints are only probabilistically enforced and since the disturbance has unbounded support, it may be that an unfortunate realization of $w_t$ drives the state far away in a region where the state constraint is strongly violated, so that no feasible control action exists to steer the state back in the region where $g(x_1, \dots,x_{M},u_0, \dots,u_{M-1}) \le 0$ holds with the required probability. Said differently, in stochastic MPC when the problem formulation \eqref{eq:prob_original} is adopted recursive feasibility cannot be guaranteed in any way, and this issue cannot be overlooked since it can be easily shown that it arises with probability one as time grows unbounded. \\

The objective of this paper is that of addressing the feasibility issue illustrated above by  introducing a suitable relaxation of problem \eqref{eq:prob_original}. This relaxation is conceived so as to adhere to the intent of the original problem formulation \eqref{eq:prob_original} as much as possible: whenever the original constraint \eqref{eq:Pconstr_y} is feasible, the original problem is maintained, while, otherwise, a new optimization problem is formulated, which is oriented towards recovering the feasibility of \eqref{eq:prob_original} as quickly as possible. As it will be shown, this is achieved by introducing an additional optimization vector $h \in \mathbb{R}^{p_y}$ that is used to relax the condition $g(x_1, \dots,x_{M},u_0, \dots,u_{M-1}) \le 0$ only for those components of the vector inequality that need to be relaxed to get feasibility, and by solving a cascade of two optimization problems with probabilistic constraints.

Admittedly, this cascade of problems can be very difficult to solve in general, since problems involving probabilistic constraints can be NP-hard. The second contribution of this paper is that of introducing a resolution scheme based on randomization in order to enhance computational tractability. Specifically, we resort to the so-called scenario approach, \cite{CC05,CC06,CG2008,CGP09}, a recently introduced randomized method that can be used to provide approximate solutions to problems with probabilistic constraints. The main advantage of the scenario approach lies in the guarantees that come attached to the returned solution and that establishes a precise link between the original problem and its approximation. In this paper, such guarantees are extended to the scenario solution to the cascade of problems discussed above, which is a non-standard setup not fully covered by the available literature (see \cite{Margellos_et_al_2015} for a contribution on cascading optimization).

With reference to the previous motivating discussion, the achievements of this paper permit one to assemble a recursive feasible MPC scheme for the stochastic system \eqref{eq:sys}, by iteratively solving the proposed relaxed optimization problem over a receding horizon.

\subsection{Paper structure}

The rest of the paper is organized as follows. Some bibliographical remarks are made in Section \ref{sec:bib_remarks}, while some compact notation is introduced in Section \ref{sec:notation}.
In Section \ref{sec:relaxation}, we formally introduce the problem relaxation, while the proposed algorithmic resolution scheme based on the scenario approach is described in Section \ref{sec:scenario}. In this section, the theoretical properties of the obtained solution are also discussed. The formal proof of these properties can be found in Section \ref{sec:proof}, while a concluding numerical example is presented in Section \ref{sec:example}.

\subsection{Bibliographical remarks}\label{sec:bib_remarks}

Alternative approaches to address control problems in presence of input and state constraints for systems affected by stochastic disturbance with unbounded support have been proposed in   \cite{Batina_et_al_2002,VanHessem2006,BB2007,Ono2008,Primbs_2007,Primbs_2009,BO2009,CACL2011}.
In \cite{Batina_et_al_2002, Primbs_2007,Primbs_2009}, state constraints are dealt with by means of a penalization term accounting for the state constraint violation so as to ensure feasibility. In \cite{VanHessem2006,Ono2008,BO2009,CACL2011}, an analytic convex relaxation of probabilistic constraints is proposed, whereas in \cite{BB2007} the problem is reformulated considering a bounded disturbance obtained by suitably cutting the tails of the disturbance distribution. 
In \cite{cannon2011stochastic,Cheng2014} stochastic uncertainty with bounded support is tackled by means of suitable probabilistic tubes.
In all these approaches, the disturbance is assumed to be a sequence of i.i.d. (independent and identically distributed) random variables. Many of them also assume that the disturbance has a Gaussian distribution, \cite{Batina_et_al_2002,VanHessem2006,BB2007,Ono2008,BO2009,CACL2011}.

This paper differs from these approaches in that a randomized-based solution is considered, which allows us to drop the independence and Gaussianity assumptions. Indeed, randomized methods have been recently developed to address design in the presence of uncertainty, making solvable problems that were otherwise deemed computationally intractable, \cite{TCD2013}. This paper differs from our previous contributions \cite{DGP2013,Franklin16} where either a term penalizing state constraint violation is added to the cost or a certain pre-defined admissible deterioration of the system performance is introduced while relaxing the state constraints.

Other approaches to constrained stochastic control for system \eqref{eq:sys} based on randomized techniques have been proposed in \cite{CAL_FAG_2011, SCFM2012, CAL_FAG_2013} under the assumption, however, that the noise has bounded support, in \cite{MSL_12} considering input constraints only, and in \cite{SFFM:2014,ZhaEtal:2014} under the assumption of recursive feasibility.

A preliminary version of this paper has been presented as a conference contribution in \cite{DGP2015}. The current submission deals with a more general setup and provides the proof of our main technical achievement.

\subsection{Compact notation}\label{sec:notation}
In order to ease the notation we define:
\begin{align*}
&\hx=
\begin{bmatrix}
x_{1} \\
x_{2}\\
\vdots \\
x_{M}
\end{bmatrix}
\quad
\hu =
\begin{bmatrix}
u_{0} \\
u_{1} \\
\vdots \\
u_{M-1}
\end{bmatrix}
 \quad
\hw =
\begin{bmatrix}
w_{0} \\
w_{1} \\
\vdots \\
w_{M-1}
\end{bmatrix}.
\end{align*}
Then, the state vector can be calculated as:
\begin{align}\label{eq:sysSR}
&\hx = \hF x_0 + \hG \hu + \hH \hw,
\end{align}
where matrices $\hF$, $\hG$ and $\hH$ are given by
\begin{align*}
&\hF = \begin{bmatrix}
A \\
A^2 \\
\vdots \\
A^M
\end{bmatrix}
 \quad
\hG = \begin{bmatrix}
B & 0_{n \times m} & \cdots & 0_{n \times m} \\
AB & B & \ddots & \vdots  \\
\vdots & \ddots & \ddots & 0_{n \times m} \\
A^{M-1}B & \cdots & AB & B
\end{bmatrix}\\
&\hH = \begin{bmatrix}
B_w & 0_{n \times n_w} & \cdots & 0_{n \times n_w} \\
AB_w & B_w & \ddots & \vdots  \\
\vdots & \ddots & \ddots & 0_{n \times n_w} \\
A^{M-1}B_w & \cdots & AB_w & B_w
\end{bmatrix}.
\end{align*}

Similarly, the disturbance feedback policy \eqref{eq:input_param} can be rewritten in the following compact form
\begin{align}\label{eq:inp_comp}
&\hu = \Gamma + \Theta \hw,
\end{align}
where we let
$$
\Gamma = \begin{bmatrix}
\gamma_0 \\
\gamma_1 \\
\vdots \\
\gamma_{M-1}
\end{bmatrix}
\Theta = \begin{bmatrix}
0_{m \times n_w} & 0_{m \times n_w} & \cdots & 0_{m \times n_w} \\
\theta_{1,0} & 0_{m \times n_w} & \ddots & \vdots  \\
\vdots & \ddots & \ddots & 0_{m \times n_w} \\
\theta_{M-1,0} & \cdots & \theta_{M-1,M-2} & 0_{m \times n_w}
\end{bmatrix}.
$$
By substituting the expression of the input in \eqref{eq:inp_comp} into \eqref{eq:sysSR}, the affine dependence of $\hx$ on the design parameters $\Gamma$ and $\Theta$ becomes clear:
\begin{align*}
&\hx = \hF x_0 + \hG\Gamma +(\hG\Theta + \hH) \hw
\end{align*}

Eventually, the nonzero components of $\Gamma$ and $\Theta$ are collected in the vector of optimization variables $\opt$, so that the following notations can be adopted:
\begin{align*}
\hu =\hu(\hw,\opt), \quad \hx =\hx(\hw,\opt),\quad J=J(\opt),
\end{align*}
which point out the dependence of input, state, and cost on the optimization vector $\opt$ and the disturbance realization $\hw$.
The constraints in \eqref{eq:Pconstr_y} then become
\begin{align*}
\mathbb{P}\{ f(\hu(\hw,\opt)) \le 0\ \wedge \ g(\hx(\hw,\opt),\hu(\hw,\opt)) \le h \} \ge 1-\epsilon.
\end{align*}

\section{Problem relaxation to ensure feasibility}\label{sec:relaxation}

In order to recover feasibility, we introduce a relaxation of the condition $g(\hx(\hw,\opt),\allowbreak \hu(\hw,\opt)) \le 0$ by substituting its right-hand side with $h \in \mathbb{R}^{p_y}$, $h$ being a new optimization variable. By doing this, the constraint involving state variables turns out to be always feasible, since it is enough to take the variable $h$ large enough. On the other hand, large values for $h$ are clearly not desired since the bigger $h$ the larger the deviation from the original constraint. To stick to the original problem formulation as much as possible $h$ should be minimized component-wise. On the other hand, one should account for the minimization of the cost function $J(\opt)$, which represents the system performance. To this purpose, the following cascade of optimization programs (two-step approach) is proposed, where $L(h)$ is an user-chosen strictly convex function of $h$, that is positive definite at $h=0$ (i.e., $L(h)> 0$, $h\ne 0$ and $L(0)=0$):
\begin{subequations} \label{eq:prob_P}
\begin{eqnarray} \label{eq:prob_Ph}
&\min_{\opt, h}\ L(h)\\
&\text{subject to: } &
\mathbb{P}\{ f(\hu(\hw,\opt)) \le 0\ \wedge \ g(\hx(\hw,\opt),\hu(\hw,\opt)) \le h \} \ge 1-\epsilon\nonumber\\
&& h \ge 0\nonumber
\end{eqnarray}
\begin{eqnarray}\label{eq:prob_PJ}
&\min_{\opt}\ J(\opt)\\
&\text{subject to: } &\mathbb{P}\left\{ f(\hu(\hw,\opt)) \le 0 \wedge g(\hx(\hw,\opt),\hu(\hw,\opt)) \le h^o\right\} \ge 1-\epsilon
\nonumber
\end{eqnarray}
\end{subequations}
where $h^o$ is the optimal value for $h$ obtained in \eqref{eq:prob_Ph}.

Problem \eqref{eq:prob_Ph} in the first step aims at determining the smallest, according to the cost $L(h)$, value of $h$ that ensures the feasibility of the probabilistic constraint 
$$
\mathbb{P}\{ f(\hu(\hw,\opt)) \le 0 \ \wedge \ g(\hx(\hw,\opt),\hu(\hw,\opt)) \le h \} \ge 1-\epsilon.
$$
A possible choice for the cost function $L(h)$ is e.g. $L(h)=h^TTh$, which allows to assign a different importance to each component of $h$ by properly choosing the positive definite matrix $T$. Note that since the cost function $L(h)$ does not depend on $\opt$, it may happen that the optimal cost $L(h^o)$ is achieved in correspondence of different choices for $\opt$, each of them leading to a possibly different value of $J(\opt)$. The second step optimization problem \eqref{eq:prob_PJ} then exploits this degree of freedom to minimize the performance cost. To this purpose,  $J(\opt)$ is minimized while the relaxed constraint $\mathbb{P}\left\{ f(\hu(\hw,\opt)) \le 0 \wedge g(\hx(\hw,\opt),\hu(\hw,\opt)) \le h^o\right\} \ge 1-\epsilon$ is enforced. Since the bound on the state condition $g(\hx(\hw,\opt),\hu(\hw,\opt))$ is fixed to $h^o$ as computed in the previous step, problem \eqref{eq:prob_PJ} does not suffer from any feasibility issue.

The cascade of problems is conceived so that, when the probabilistic constraint in \eqref{eq:prob_original} is infeasible, the control action is basically designed according to \eqref{eq:prob_Ph} so as to recover the feasibility as quickly as possible (minimization of $L(h)$). In this case, \eqref{eq:prob_PJ} provides just a refinement of the solution. The addition of the constraint $h \geq 0$ in \eqref{eq:prob_Ph} ensures that the constraint relaxation in \eqref{eq:prob_PJ}, component by component, cannot become tighter than the original constraint in \eqref{eq:prob_original}, and for those components not requiring any relaxation \eqref{eq:prob_PJ} pursues the goal of minimizing $J(\opt)$ as in \eqref{eq:prob_original}. In particular, whenever \eqref{eq:prob_original} is already feasible, program \eqref{eq:prob_Ph} simply returns $h^o = 0$ and the original problem \eqref{eq:prob_original} is recovered in \eqref{eq:prob_PJ}.

Overall, the cascade of problems \eqref{eq:prob_P} returns a solution given by the pair $(\opt^o,h^o)$, where $\opt^o$ determines the control action to be implemented and $h^o$ is the probabilistically guaranteed bound for the state constraint. Note that the value $h^o$ computed in the first step optimization problem can be inspected to evaluate the mismatch with respect to the original state constraint.

\section{Scenario-based resolution scheme}\label{sec:scenario}

As problems \eqref{eq:prob_Ph} and \eqref{eq:prob_PJ} are, in general, hard to solve because of the presence of a probabilistic constraint, we propose to tackle them by means of a randomized scheme which is in the vein of the so-called scenario approach, \cite{CC05,CC06,CG2008,CGP09}. The proposed scheme allows to recover computational tractability at the price of introducing some approximation. However, by exploiting the scenario approach, which is here extended to the cascade of problems \eqref{eq:prob_Ph} and \eqref{eq:prob_PJ}, precise probabilistic guarantees on the feasibility of the achieved solution are also provided.

The idea of the scenario approach is to consider $N$ disturbance realizations, each extracted according to the underlying probability distribution $\mathbb P$:
\begin{align*}
&\mathbf{w}^{(1)}=\left[w_0^{(1)}\ w_1^{(1)}\ \ldots\ w_{M-1}^{(1)}\right]\\
&\mathbf{w}^{(2)}=\left[w_0^{(2)}\ w_1^{(2)}\ \ldots\ w_{M-1}^{(2)}\right]\\
&\ \vdots\\
&\mathbf{w}^{(N)}=\left[w_0^{(N)}\ w_1^{(N)}\ \ldots\ w_{M-1}^{(N)}\right].
\end{align*}
Then, the probabilistic constraint in \eqref{eq:prob_Ph} and \eqref{eq:prob_PJ} are replaced by $N$ non-probabilistic constraints, one for each disturbance realizations. More precisely, we have the following cascade of problems that can be seen as a randomized counterpart of the cascade of problems \eqref{eq:prob_P}:
\begin{subequations} \label{eq:prob_Scenario}
\begin{eqnarray} \label{eq:prob_Sh}
&\min_{\lambda, h} L(h) \\
&\text{subject to: } &f(\hu(\hw^{(k)},\lambda)) \le 0, \nonumber \\
& &g(\hx(\hw^{(k)},\lambda)),\hu(\hw^{(k)},\lambda)) \le h, \quad  k=1, \ldots, N, \nonumber \\
& &h \ge 0 . \nonumber
\end{eqnarray}
\begin{eqnarray} \label{eq:prob_SJ}
& \min_{\lambda}   J(\lambda)\\
&\text{subject to: }  & f(\hu(\hw^{(k)},\lambda)) \le 0, \nonumber \\
& & g(\hx(\hw^{(k)},\lambda)),\hu(\hw^{(k)},\lambda)) \le h^\star, \quad  k=1, \ldots, N, \nonumber
\end{eqnarray}
\end{subequations}
where $h^\star$ is the optimal value of $h$ obtained in \eqref{eq:prob_Sh}.

Problems \eqref{eq:prob_Sh} and \eqref{eq:prob_SJ} are convex and have a finite number of constraints, hence they can be efficiently solved by means of standard solvers. Note that as the constraints are convex and the cost function $L(h)$ is assumed to be strictly convex with respect to its argument, problem \eqref{eq:prob_Sh} uniquely determines the value of $h^\star$, similarly, thanks to the strict convexity of $J(\opt)$, the solution of problem \eqref{eq:prob_SJ}, say $\lambda^\star$, is unique.

The same interpretation we had for the cascade of problems \eqref{eq:prob_P} in Section \ref{sec:relaxation} applies to the cascade of problems \eqref{eq:prob_Scenario}: indeed, the solution of \eqref{eq:prob_Scenario} defined by the pair $(\lambda^\star,h^\star)$ is the empirical counterpart of the solution of \eqref{eq:prob_P}. It is worth noticing that, as the pair $(\opt^\star,h^\star)$ is feasible and optimal for \eqref{eq:prob_Sh}, the second step optimization problem \eqref{eq:prob_SJ} can be regarded as a tie break rule by means of which the solution that minimizes the cost $J(\opt)$ is chosen among the possible multiple solutions in $\lambda$ of the first step optimization problem \eqref{eq:prob_Sh}.

\subsection{Probabilistic constraint feasibility of the scenario solution}

Resorting to \eqref{eq:prob_Scenario} allows to enhance computational tractability, but a question arises on whether the obtained approximated solution $(\lambda^\star,h^\star)$ is satisfactory or not as far as the original probabilistically constrained cascade of problems is concerned. In particular, we are interested in studying the feasibility of the obtained scenario-based solution for the probabilistic constraint
\begin{equation} \label{eq:prob_constr}
\mathbb{P}\left\{ f(\hu(\hw,\opt)) \le 0 \wedge g(\hx(\hw,\opt),\hu(\hw,\opt)) \le h \right\} \ge 1-\epsilon,
\end{equation}
so as to provide a connection between $(\lambda^\star,h^\star)$ and the original cascade of problems \eqref{eq:prob_P}.

This question pertains to the theory of the scenario approach, which provides in a number of different setups guarantees on the feasibility of the scenario solution for the original probabilistic constraint as long as $N$ is suitably chosen, \cite{CC05,CC06,CG2008,CG2011,GC2013}. The best available result is that of \cite{CG2008} which, however, does not directly apply to the cascade of problems \eqref{eq:prob_Scenario}. Indeed, the result in  \cite{CG2008} are proven for scenario optimization problems whose solution is determined by a specific tie break rule, \cite[Section 2.1 point 5]{CG2008}, which may not correspond to the one determined by the cascade of problems \eqref{eq:prob_Scenario}.

The results on cascading optimization in \cite{Margellos_et_al_2015} apply to this context but the resulting bound on $N$ is more conservative that the one
in the following theorem, which provides an extension of the result in \cite{CG2008} to the current framework, and is tailored to the specific setting at hand, where the second problem in the cascade can be interpreted as a tie breaking rule for the first one.

\begin{theorem} \label{th:main-result}
Let $\beta \in (0,1)$ be a user-chosen confidence parameter. If the number of extracted disturbance realizations $N$ is chosen so as to satisfy
\begin{align} \label{eq:betaN}
\sum_{i=0}^{d-1}\binom{N}{i}\epsilon^i(1-\epsilon)^{N-i} \leq \beta,
\end{align}
where $d$ is the dimensionality of $(\lambda,h)$, then it holds with confidence at least $1-\beta$ that
\begin{align*}
\mathbb{P} \left\{ f(\hu(\hw,\lambda^\star)) \le 0 \ \wedge \ g(\hx(\hw,\lambda^\star),\hu(\hw,\lambda^\star)) \le h^\star \right\} \geq 1-\epsilon. \\
\end{align*}
\end{theorem}

\emph{Proof:} see next Section \ref{sec:proof}. \\

The theorem states that with high confidence $1-\beta$ the solution $(\lambda^\star,h^\star)$ achieved solving the scenario cascade of problems \eqref{eq:prob_Scenario} is feasible for the original probabilistic constraint \eqref{eq:prob_constr} in \eqref{eq:prob_P}. Note that the presence of the confidence parameter $\beta$ is intrinsic and is related to fact that the obtained solution depends on the random extraction $\hw^{(1)},\ldots,\hw^{(N)}$: $\beta$ is needed to account for the possibility that a not representative enough sample $\hw^{(1)},\ldots,\hw^{(N)}$ is seen. However exploiting the results in \cite{AlamoTempo15} it can be shown that the number of required samples $N$ according to \eqref{eq:betaN} scales logarithmically with $1/\beta$. Hence $\beta$ can be chosen to be very small such as $10^{-5}$ or $10^{-7}$ without affecting $N$ too much, so that the fact that the achieved solution $(\lambda^\star,h^\star)$ satisfies the probabilistic constraint \eqref{eq:prob_constr} in \eqref{eq:prob_P} can be taken for granted.

\section{Proof of Theorem \ref{th:main-result}} \label{sec:proof}

This proof is inspired by the proof of  point 5  in \cite[Section 2.1]{CG2008}, which was not included in \cite{CG2008} but provided by the authors of \cite{CG2008} on request.

For a given $(\lambda,h)$, define the violation probability of $(\lambda,h)$ as
\begin{align*}
V(\lambda,h) := \mathbb{P} \Big\{
f(\hu(\hw,\lambda)) > 0 \; \vee \; g(\hx(\hw,\lambda),\hu(\hw,\lambda)) > h
\Big\}& \\
= 1 - \mathbb{P} \Big\{
f(\hu(\hw,\lambda)) \leq 0 \ \wedge \ g(\hx(\hw,\lambda),\hu(\hw,\lambda)) \leq h
\Big\},
\end{align*}
where $\vee$ stands for ``or''.
Then, Theorem \ref{th:main-result} amounts to showing that
\begin{equation} \label{eq:P^N(V>eps)<=beta}
\mathbb{P}^N \{ V(\lambda^\star,h^\star) > \epsilon \} \leq \beta,
\end{equation}
where $\mathbb{P}^N$ is the product probability underlying the independent extraction of the sample $\hw^{(1)},\ldots,\hw^{(N)}$ based on which the solution $(\lambda^\star,h^\star)$ is computed.

Consider the following auxiliary scenario program
\begin{eqnarray} \label{eq:prob_aux}
&\min_{\lambda, h} L(h) + \frac{1}{n} J(\lambda) \\
&\text{subject to:}  & f(\hu(\hw^{(k)},\lambda)) \le 0, \nonumber \\
& & g(\hx(\hw^{(k)},\lambda),\hu(\hw^{(k)},\lambda)) \le h, \quad  k=1 \ldots N, \nonumber \\
& & h \ge 0, \nonumber
\end{eqnarray}
for $n=1,2,\ldots$, and denote by $(\lambda^\star_n,h^\star_n)$ its optimal solution, which exists and is unique, since: i. the cost function $L(h) + \frac{1}{n} J(\lambda)$ has compact level sets for every $n \geq 1$ thanks to the strict convexity of $L$ and on $J$; ii. the optimization feasibility domain defined by the constraints in \eqref{eq:prob_aux} is close and nonempty. The following two properties hold.
\begin{itemize}
\item[1.] For every $n \geq 1$, it holds that
\begin{equation} \label{eq:aux_scen_ok}
\mathbb{P}^N \{ V(\lambda^\star_n,h^\star_n) > \epsilon \} \leq \beta.
\end{equation}
\item[2.] For every multisample $\hw^{(1)},\ldots,\hw^{(N)}$, the solution to \eqref{eq:prob_aux} converges to the solution to \eqref{eq:prob_Scenario}, namely,
\begin{equation} \label{eq:convergence}
(\lambda^\star_n,h^\star_n) \to (\lambda^\star,h^\star) \text{ as } n \to \infty.
\end{equation}
\end{itemize}

\noindent {\it Proof of Property 1} \\

By adding a slack variable $v \in \mathbb R$, problem \eqref{eq:prob_aux} can be rewritten in epigraphic form as:
\begin{eqnarray} \label{eq:prob_aux_epigraphic}
&\min_{\lambda, h, v}  v \\
& \text{subject to:} & f(\hu(\hw^{(k)},\lambda)) \le 0, \nonumber \\
& & g(\hx(\hw^{(k)},\lambda),\hu(\hw^{(k)},\lambda)) \le h, \quad  k=1 \ldots N, \nonumber \\
& & h \ge 0, \nonumber \\
& & L(h) + \frac{1}{n} J(\lambda) \leq v. \nonumber
\end{eqnarray}
The solution to problem \eqref{eq:prob_aux_epigraphic} is still unique, and the assumptions of Theorem 2.4 in \cite{CG2008} are satisfied. An application of this theorem gives
$$
\mathbb{P}^N \{ V(\lambda^\star_n,h^\star_n) > \epsilon \} \leq \sum_{i=0}^{d}\binom{N}{i}\epsilon^i(1-\epsilon)^{N-i},
$$
where we have $d$ in place of $d-1$ because in \eqref{eq:prob_aux_epigraphic} the number of optimization variables has been augmented by $1$ and is equal to $d+1$. On the other hand, since the slack variable $v$ does not enter the expression defining it, the constraint
$$
\{ \lambda,h,v : \ f(\hu(\hw,\lambda)) \le 0 \; \wedge \; g(\hx(\hw,\lambda),\hu(\hw,\lambda)) \le h \}
$$
is, irrespective of $\hw$, a cylindroid infinitely extended along the $v$ direction. This entails that the family (with respect to the variability of $\hw$) of constraints above has a so-called support rank equal to $d$, according to Definition 3.6 of \cite{SfagM2013} (see also \cite{ZhaEtal:2015}). The conclusion that
$$
\mathbb{P}^N \{ V(\lambda^\star_n,h^\star_n) > \epsilon \} \leq \sum_{i=0}^{d-1}\binom{N}{i}\epsilon^i(1-\epsilon)^{N-i}
$$
then follows by invoking the observation made in \cite{SfagM2013} that Theorem 2.4 of \cite{CG2008} still applies by replacing the optimization domain dimensionality with the support rank (see Lemma 3.8). \qed \\

\noindent {\it Proof of Property 2} \\

To show that $(\lambda^\star_n,h^\star_n) \to (\lambda^\star,h^\star)$ as $n \to \infty$, consider the sets
\begin{align*}
\mathcal{H}_n = \Big\{ (\lambda,h): \  (\lambda,h) \text{ is feasible for \eqref{eq:prob_aux} and}  L(h) + \frac{1}{n} J(\lambda) \leq L(h^\star) + \frac{1}{n} J(\lambda^\star) \Big\},
\end{align*}
for $n=1,2,\ldots$. In words, $n$ by $n$, $\mathcal{H}_n$ is the set of all feasible points for \eqref{eq:prob_aux} that also belong to the smallest level set of the cost function of \eqref{eq:prob_aux} containing the solution $(\lambda^\star,h^\star)$ of \eqref{eq:prob_Scenario}. Note that, while the level set changes with $n$, the feasibility domain of \eqref{eq:prob_aux} remains the same for all $n$ and it coincides with the feasibility domain of \eqref{eq:prob_Sh}. This entails that $(\lambda^\star,h^\star)$ belongs to $\mathcal{H}_n$ for all $n$, showing also that $\mathcal{H}_n$ is nonempty. Moreover, $n$ by $n$, we have that
\begin{equation} \label{eq:sol_inclusion}
(\lambda^\star_n,h^\star_n) \in \mathcal{H}_n,
\end{equation}
because $(\lambda^\star_n,h^\star_n)$ is feasible for \eqref{eq:prob_aux}, and, being also optimal, its cost value must be better than that of  $(\lambda^\star,h^\star)$, which is the second condition defining $\mathcal{H}_n$.

\begin{figure}
\centering
\vspace{0.2cm}
\includegraphics[trim=4.4cm 9cm 2cm 7cm,clip=true,scale=0.6]{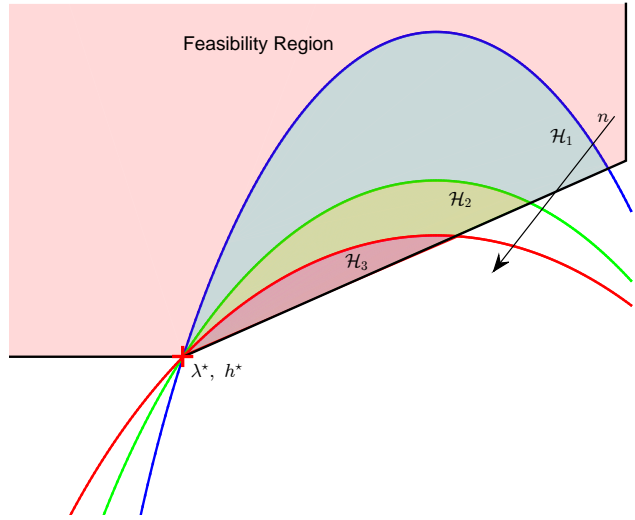}
\caption{The sets $\mathcal H_n$'s in a simple case ($h,\opt \in \mathbb{R}$, $L(h)=h^2$, $J(\opt)=3\opt^2$).}\label{fig:Proof}
\end{figure}

A fundamental property of the family of sets $\mathcal{H}_n$ is that
\begin{equation} \label{eq:H_n_inclusion}
\mathcal{H}_{1} \supseteq \mathcal{H}_{2} \supseteq \cdots \supseteq \mathcal{H}_{n} \supseteq \mathcal{H}_{n+1} \supseteq \cdots,
\end{equation}
as pictorially depicted in Fig. \ref{fig:Proof}. To show \eqref{eq:H_n_inclusion}, suppose that a $(\bar \lambda,\bar h)$ belongs to $\mathcal{H}_{n+1}$. From
$$
L(\bar{h}) + \frac{1}{n+1} J(\bar{\lambda}) \leq L(h^\star) + \frac{1}{n+1} J(\lambda^\star)
$$
it follows that $J(\bar{\lambda}) \leq (n+1)( L( h^\star) - L(\bar h)) + J(\lambda^\star)$. Whence,
\begin{align*}
{ L(\bar{h}) + \frac{1}{n} J(\bar{\lambda})}
& \leq   L(\bar{h}) + \frac{n+1}{n}(L( h^\star) - L(\bar{h})) +  \frac{1}{n} J(\lambda^\star) \\
& =  L(h^\star) + \frac{1}{n} (L(h^\star) - L(\bar{h})) +  \frac{1}{n} J(\lambda^\star) \\
& \leq  L( h^\star) + \frac{1}{n} J(\lambda^\star),
\end{align*}
where the last inequality follows because $L(h^\star) - L(\bar{h}) \leq 0$ being $L(h^\star)$ the lowest among feasible points by the definition of $h^\star$.
This shows that $(\bar \lambda,\bar h) \in \mathcal{H}_{n}$ too, that is, \eqref{eq:H_n_inclusion} holds.

From \eqref{eq:sol_inclusion} and \eqref{eq:H_n_inclusion}, we have that $(\lambda^\star_n,h^\star_n) \in \mathcal{H}_1$, $\forall n$. Set $\mathcal{H}_1$ is compact, being the intersection of the feasibility domain of \eqref{eq:prob_Sh}, which is close, with a level set of $L(h) + \frac{1}{n} J(\lambda)$, which is compact thanks to the assumptions of strict of convexity of $L$ and $J$. It then follows that the sequence $(\lambda^\star_n,h^\star_n)$ have limit points, which are feasible for \eqref{eq:prob_Sh}. For simplicity, assume that there is just one, say $(\lambda^\star_\infty,h^\star_\infty)$, so that the sequence $(\lambda^\star_n,h^\star_n)$ is convergent to $(\lambda^\star_\infty,h^\star_\infty)$. If not, simply repeat the argument that follows to each limit point and the corresponding convergent subsequence.

From \eqref{eq:sol_inclusion} and the definition of $\mathcal{H}_n$, we have that
$$
L(h^\star_n) \leq L(h^\star) + \frac{1}{n} \left[ J(\lambda^\star) - J(\lambda^\star_n) \right],
$$
which in turn implies that
\begin{eqnarray*}
L(h^\star_\infty) & = & \lim_{n \to \infty} L(h^\star_n) \\
& \leq & L(h^\star) + \lim_{n \to \infty} \frac{1}{n} \left[ J(\lambda^\star) - J(\lambda^\star_n) \right] \\
& = & L(h^\star).
\end{eqnarray*}
Yet, being $L(h^\star)$ minimal, it cannot be that a strict inequality holds, so that eventually $L(h^\star_\infty) = L(h^\star)$. If $h^\star_\infty \neq h^\star$, then $(\frac{1}{2}\lambda^\star+\frac{1}{2}\lambda^\star_\infty,\frac{1}{2}h^\star+\frac{1}{2}h^\star_\infty)$ would be feasible for \eqref{eq:prob_Sh} thanks to the convexity of the feasible domain, while the strict convexity of $L(h)$ would give
\begin{eqnarray*}
L\left( \frac{1}{2}h^\star+\frac{1}{2}h^\star_\infty \right) & < & \frac{1}{2} L(h^\star) + \frac{1}{2} L(h^\star_\infty) \\
& = & L(h^\star),
\end{eqnarray*}
so contradicting the minimality of $L(h^\star)$. Hence, $h^\star_\infty = h^\star$.

From $(\lambda^\star_n,h^\star_n) \in \mathcal{H}_1$, we have that
$J(\lambda^\star_n) \leq L(h^\star) - L(h^\star_n) +  J(\lambda^\star)$ which, taking the limit, gives
\begin{eqnarray*}
J(\lambda^\star_\infty) & \leq & \lim_{n \to \infty} L(h^\star) - L(h^\star_n) +  J(\lambda^\star) \\
& = & J(\lambda^\star).
\end{eqnarray*}
Plainly, it must be that $J(\lambda^\star_\infty) = J(\lambda^\star)$, for, otherwise, being $\lambda^\star_\infty$ feasible for \eqref{eq:prob_SJ}, $J(\lambda^\star_\infty) < J(\lambda^\star)$ would contradict the minimality of $J(\lambda^\star)$. Moreover, if $\lambda^\star_\infty \neq \lambda^\star$, then $\frac{1}{2}\lambda^\star+\frac{1}{2}\lambda^\star_\infty$ would be feasible for \eqref{eq:prob_SJ}, and, because of the strict convexity of $J(\lambda)$ we would have
$$
J(\frac{1}{2}\lambda^\star+\frac{1}{2}\lambda^\star_\infty) < \frac{1}{2}J(\lambda^\star) + \frac{1}{2}J(\lambda^\star_\infty) = J(\lambda^\star),
$$
contradicting again the minimality of $J(\lambda^\star)$. Hence, $\lambda^\star_\infty = \lambda^\star$, and this concludes the proof of Property 2. \qed \\

We want now to capitalize on \eqref{eq:aux_scen_ok} and \eqref{eq:convergence} to show that \eqref{eq:P^N(V>eps)<=beta} holds. To this purpose, start by fixing a sample $\hw^{(1)},\ldots,\hw^{(N)}$ such that $V(\lambda^\star,h^\star) > \epsilon$, which, we recall, means that
$$
\mathbb{P} \Big\{
f(\hu(\hw,\lambda^\star)) > 0 \; \vee \; g(\hx(\hw,\lambda^\star),\hu(\hw,\lambda^\star)) > h^\star
\Big\} > \epsilon.
$$
By continuity of $f$ and $g$, this implies that
$$
\mathbb{P} \Big\{
f(\hu(\hw,\lambda)) > 0 \; \vee \; g(\hx(\hw,\lambda),\hu(\hw,\lambda)) > h^\star \Big\} > \epsilon,
$$
for all $(\lambda,h): \| (\lambda,h) - (\lambda^\star,h^\star) \| \leq r$ for a radius $r$ small enough, and, since  $(\lambda^\star_n,h^\star_n) \to (\lambda^\star,h^\star)$ so that $\| (\lambda,h) - (\lambda^\star,h^\star) \| \leq r$ for all $n$ bigger than a suitable $\bar n$, we can conclude that
\begin{align} \label{eq:V(^star_n)>eps}
V(\lambda^\star_n,h^\star_n) = \mathbb{P} \Big\{
f(\hu(\hw,\lambda^\star_n)) > 0 \ \vee \ g(\hu(\hw,\lambda^\star_n),\hx(\hw,\lambda^\star_n)) > h^\star_n \Big\} > \epsilon,
\end{align}
for $n > \bar n$. If we now let $\hw^{(1)},\ldots,\hw^{(N)}$ vary and we consider the indicator function $\mathbb{I}_{\{ \hw^{(1)},\ldots,\hw^{(N)}: \; V(\lambda^\star_n,h^\star_n) > \epsilon\}}$, then \eqref{eq:V(^star_n)>eps} yields
$$
\mathbb{I}_{\{ V(\lambda^\star,h^\star) > \epsilon\}} \cdot \mathbb{I}_{\{ V(\lambda^\star_n,h^\star_n) > \epsilon\}} \xrightarrow[n \to \infty]{} \mathbb{I}_{\{ V(\lambda^\star,h^\star) > \epsilon\}},
$$
for all possible realizations of $\hw^{(1)},\ldots,\hw^{(N)}$. Applying the Lebesgue dominated convergence theorem gives
\begin{align*}
{ \lim_{n \to \infty} \mathbb{P}^N \{ V(\lambda^\star_n,h^\star_n) > \epsilon \} }
& =  \lim_{n \to \infty} \int \mathbb{I}_{ \{ V(\lambda^\star_n,h^\star_n) > \epsilon \} } \mathbb{P}^N \{ \dd \hw^{(1)},\ldots, \dd \hw^{(N)} \} \\
& \geq  \lim_{n \to \infty} \int \mathbb{I}_{ \{ V(\lambda^\star,h^\star) > \epsilon\} } \cdot \mathbb{I}_{ \{ V(\lambda^\star_n,h^\star_n) > \epsilon \} } \mathbb{P}^N \{ \dd \hw^{(1)},\ldots, \dd \hw^{(N)} \} \\
& =  \int \mathbb{I}_{ \{ V(\lambda^\star,h^\star) > \epsilon\} } \mathbb{P}^N \{ \dd \hw^{(1)},\ldots, \dd \hw^{(N)} \} \\
& =  \mathbb{P}^N \{ V(\lambda^\star,h^\star) > \epsilon \}.
\end{align*}
Hence, $\mathbb{P}^N \{ V(\lambda^\star,h^\star) > \epsilon \} \leq \lim_{n \to \infty} \mathbb{P}^N \{ V(\lambda^\star_n,h^\star_n) > \epsilon \}$, and since $\mathbb{P}^N \{ V(\lambda^\star_n,h^\star_n) > \epsilon \} \leq \beta$, $\forall n$, \eqref{eq:P^N(V>eps)<=beta} remains proven. \qed

\section{Numerical Example}\label{sec:example}

In this section we apply the proposed approach to a numerical example inspired by \cite{CACL2011}.

\begin{figure}[!h]
\centering
\includegraphics[scale=0.7]{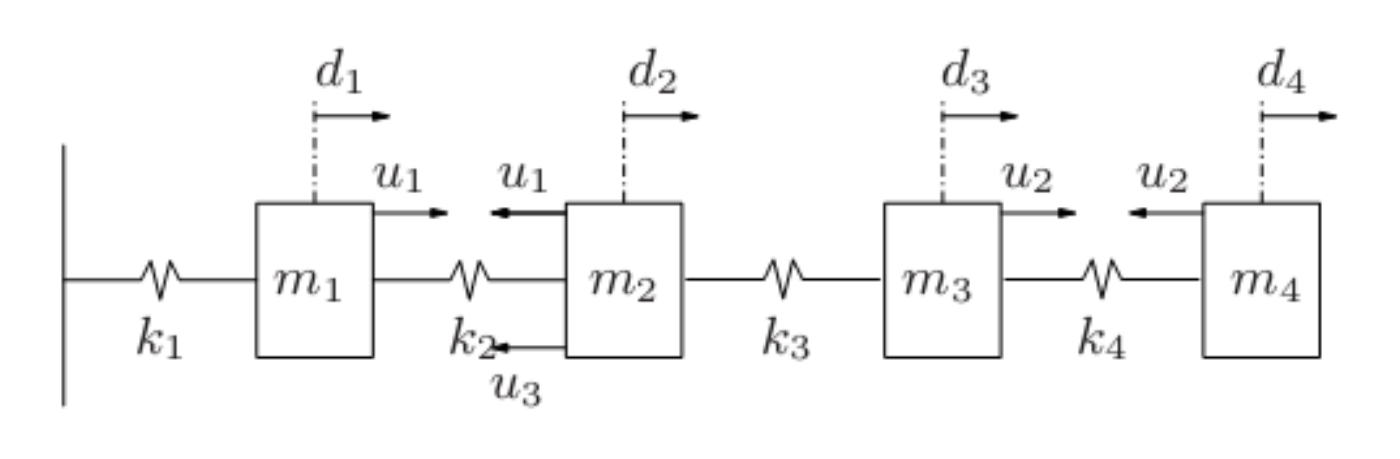}
\caption{Scheme of the mechanical system.}
\label{fig:mech_sys}
\end{figure}

The mechanical system composed by 4 masses and 4 springs depicted in Fig. \ref{fig:mech_sys} is considered.
Masses and stiffness coefficients of springs are all equal to 1. The state of the system is formed by the displacements of masses with respect to nominal positions and their derivatives, that is, $x=[d_1,\ d_2,\ d_3,\ d_4,\ \dot d_1,\ \dot d_2,\ \dot d_3,\ \dot d_4]^T$. The control inputs $u_1$, $u_2$, $u_3$ are instead the forces acting on the masses shown in Fig. \ref{fig:mech_sys}. The continuous-time system equations are easily derived. The system dynamics is then discretized assuming that the input is kept constant in the interval $[t,\ t+T_s)$, with $T_s=1$ s, so obtaining a system as in \eqref{eq:sys}. A stochastic additive disturbance affecting both the masses displacements and speeds is supposed to be also present, resulting after discretization in $w\sim WGN(0,\ I_4)$ and $B_w=[0.5I_4\ I_4]^T$. The initial condition of the system is $x_0=[10,\ -10,\ 10,\ -10,\ 0,\ 0,\ 0,\ 0]^T$.

The goal is to keep the masses close to their nominal position counteracting the action of the disturbance. We do not consider constraints on the input, but we enforce a state constraint requiring that the maximum speed of each mass keeps below a given bound.

To this purpose, we choose an average quadratic cost function as in \eqref{eq:cost} over a finite horizon $M=8$, where the matrices $Q$ and $R$ are set so as to penalize deviations from the nominal positions:
\begin{align*}
Q=\begin{bmatrix}
I_4 & 0_{4\times4} \\
0_{4\times4} & 0_{4\times4}
\end{bmatrix} \quad
R=10^{-6}I_3,
\end{align*}
The constraints on the speed of the masses are instead formulated as follows:
\begin{align} \label{eq:ex_stateconstr}
\|Cx_i \|_{\infty} \le 10 \qquad i=1,\ldots, M,
\end{align}
where $C=[0_{4\times4}\ I_4]$. To deal with the presence of the disturbance $w$, the constraint is enforced in probability with $\epsilon=0.1$.

Following the approach of Section \ref{sec:relaxation}, the optimization variables $h_i$, $i=1,\ldots, M$, are introduced so as to ensure feasibility of the probabilistic constraint. We set $\beta=10^{-6}$ resulting in $N=4614$ according to  \eqref{eq:betaN}. Eventually, the cascade of problems \eqref{eq:prob_Scenario} is solved with $L(h)=h^Th$.

Numerical results show that the bound on the state constraint cannot be enforced for the first 2 time steps. Indeed solving problem \eqref{eq:prob_Sh} gave that the smallest bound preserving feasibility is $h_1^\star=1.62$, $h_2^\star=1.08$, while for the other time steps we had $h_i^\star=0$, $i=3,\ldots, M$. The cost $J(\opt^\star)$ achieved solving problem \eqref{eq:prob_SJ} was $2305.55$. Some Monte-Carlo simulations revealed that the probabilistic constraint \eqref{eq:prob_constr} was satisfied by the achieved solution $(\opt^\star,h^\star)$ as it was expected given Theorem \ref{th:main-result}.

In order to better evaluate the performance of the obtained scenario control policy, we compared it against a finite horizon LQ controller, which was designed according to the following cost function:
\begin{align*}
&J_{LQ}=\mathbb{E}\left[\sum_{t=1}^{M} x_t^TQ_{LQ}x_t + \sum_{t=0}^{M-1} u_t^TR_{LQ}u_t \right],\\
&Q_{LQ}=\begin{bmatrix}
q_JI_4 & 0_{4\times4}\\
0_{4\times4} & q_{L}I_4
\end{bmatrix}\quad
R_{LQ}=10^{-6}I_3,
\end{align*}
where the weights $q_{J}$ and $q_{L}$ are degrees of freedom to tune the relative importance between displacements and speeds. In this way the LQ controller can partially account for the requirement on the masses speed.

\begin{figure}[th!]
\subfloat[Scenario-based solution, $10+h^\star_1$ blue dash-dotted line, $10+h^\star_2$ red dash-dotted line.]{\includegraphics[trim=1.5cm 8.1cm 2cm 8.6cm, clip=true,scale=0.33]{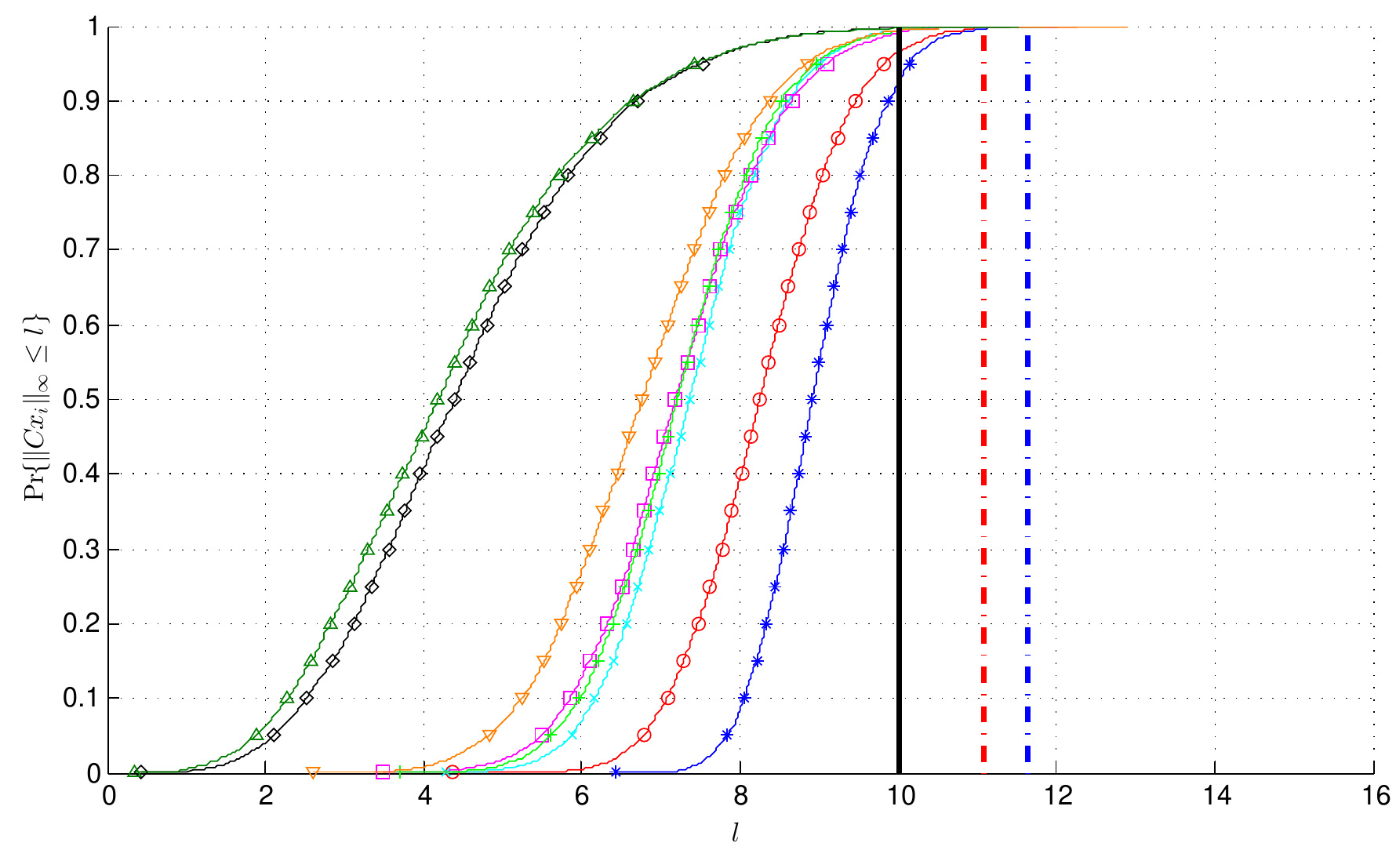}}
\subfloat[LQ controller $q_J=1$, $q_{L}=0$]{\includegraphics[trim=1.6cm 8.1cm 2cm 8.6cm, clip=true,scale=0.33]{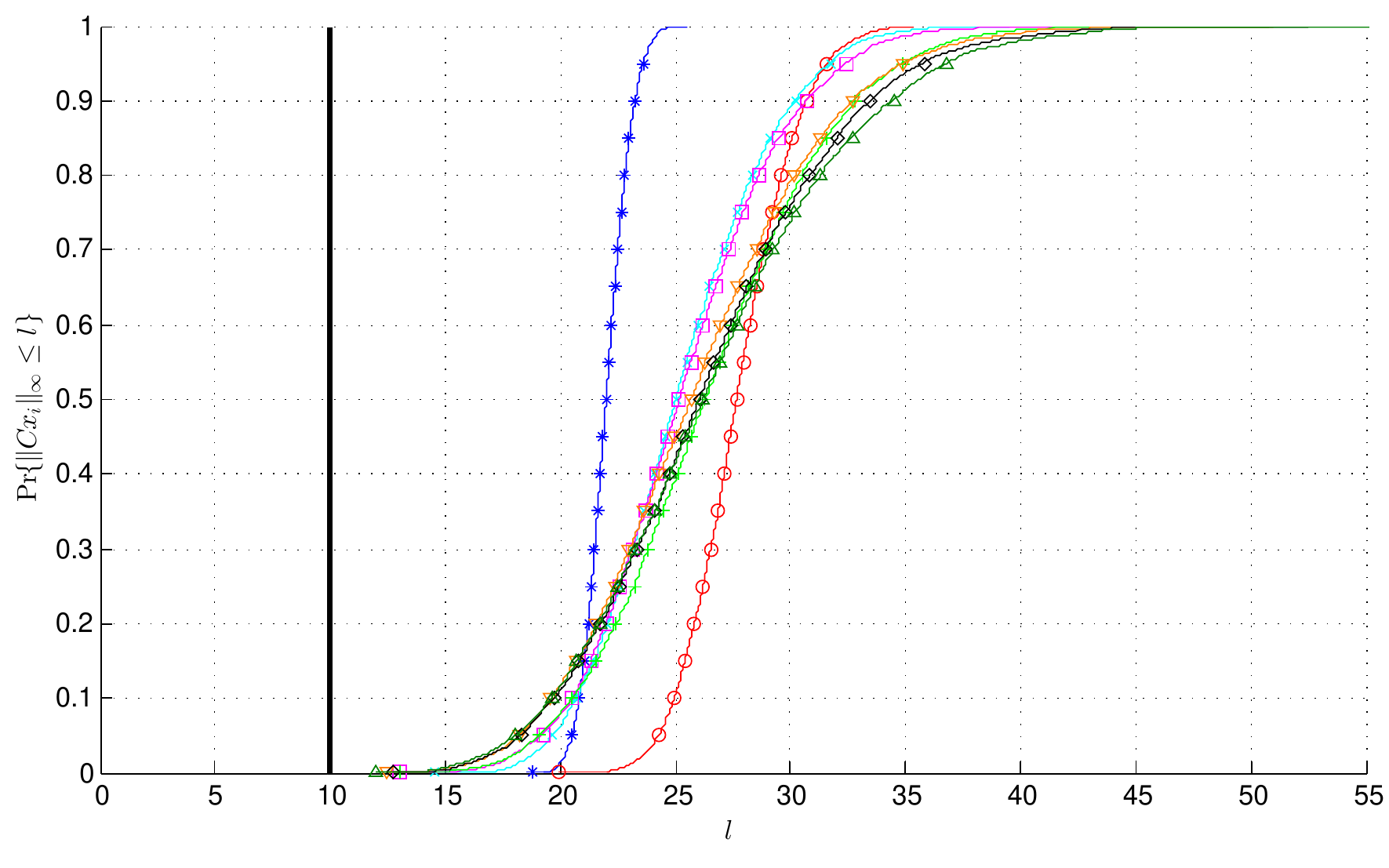}}
\subfloat{\includegraphics[
clip=true,scale=.4]{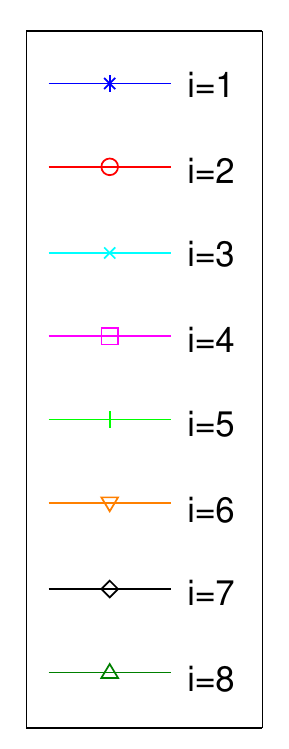}} \\
\subfloat[LQ controller $q_J=0$, $q_{L}=1$]{\includegraphics[trim=1.6cm 8.1cm 2cm 8.6cm, clip=true,scale=0.33]{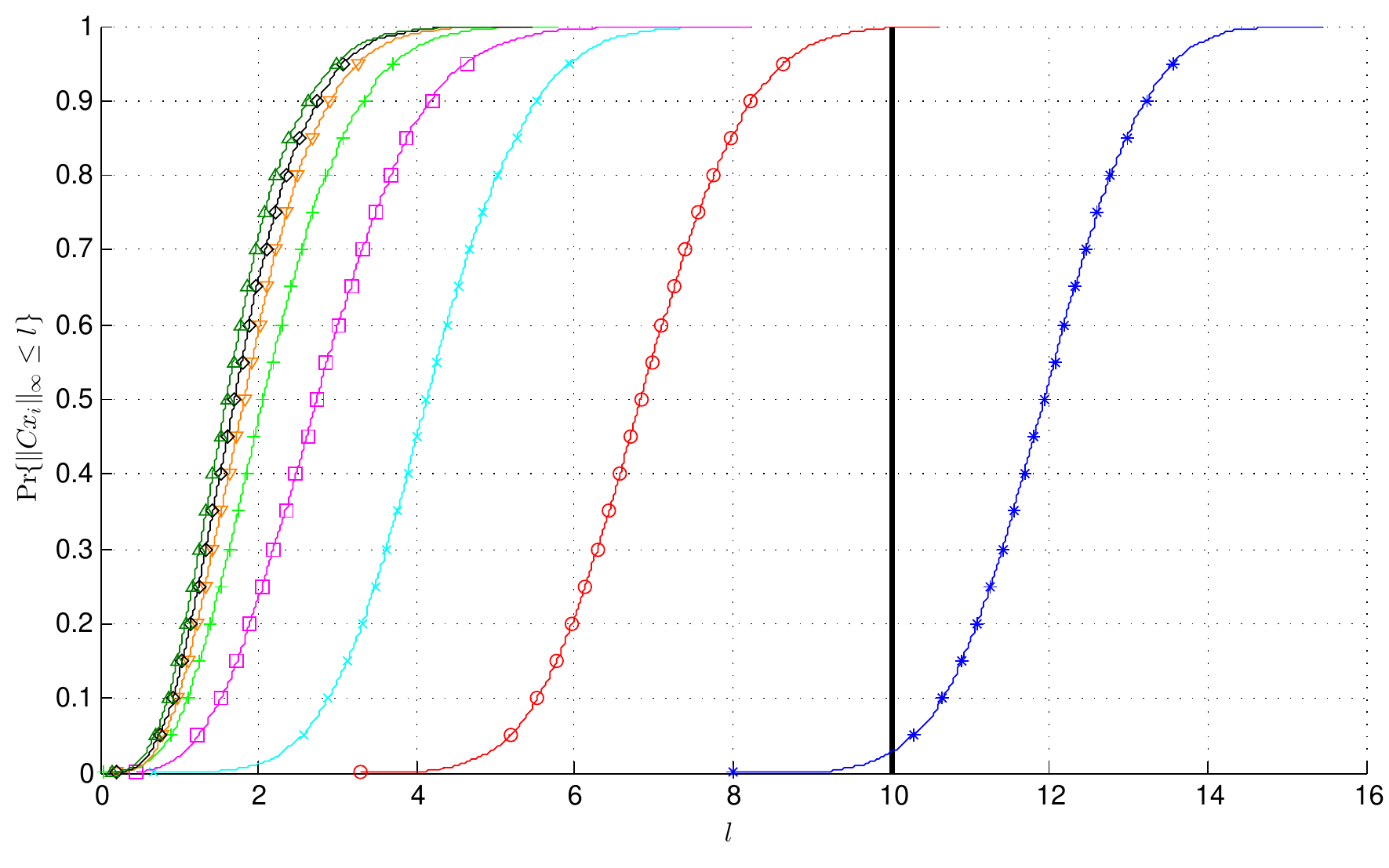}}
\subfloat[LQ controller $q_J=0.2$, $q_{L}=9$]{\includegraphics[trim=1.6cm 8.1cm 2cm 8.6cm, clip=true,scale=0.33]{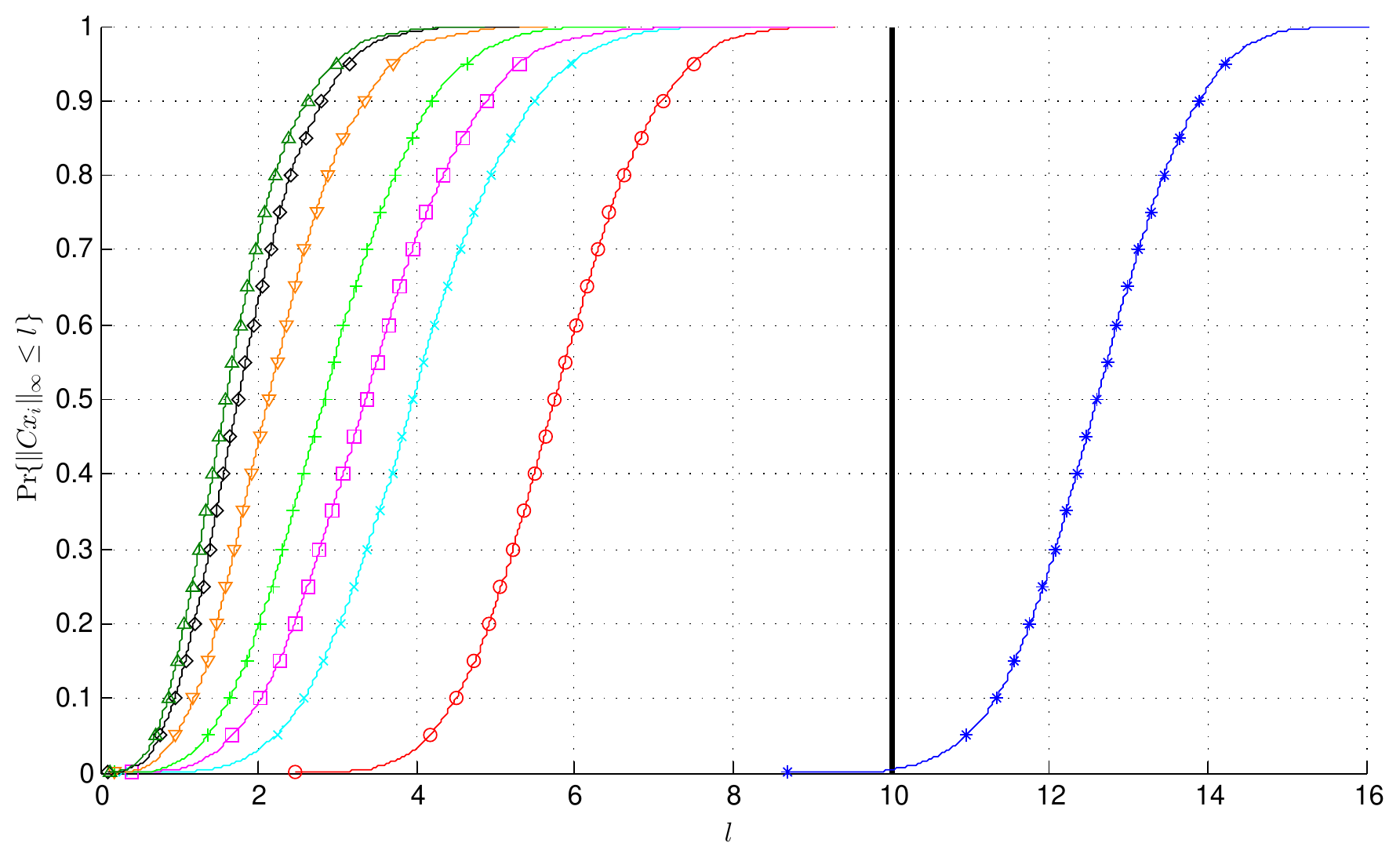}}
\caption{Cumulative probability distributions of $\|Cx_i\|_{\infty}$, $i=1,\ldots, 8$ for the scenario-based solution and for the LQ controllers; the feasible domain of the speed constraint is represented by the left side of the marked vertical solid line.}
\label{fig:distr}
\end{figure}

The comparison of the performance of the scenario-based control policy and that of the LQ controller for different choices of $q_{J}$ and $q_{L}$ is displayed in Table \ref{tab:exaPERF} and in Fig. \ref{fig:distr}.
Specifically, Table \ref{tab:exaPERF} reports the achieved cost $J$ and the actual probability of violation $\tilde \epsilon$ of the original state constraint \eqref{eq:ex_stateconstr} (computed via Monte Carlo simulations), while Fig. \ref{fig:distr} depicts the cumulative probability distributions of $\|Cx_i\|_{\infty}$, $i=1,\ldots,M$, (again, as obtained by means of Monte Carlo simulations).

\begin{table}[!h]
  \centering
  \caption{}
    \begin{tabular}{|cc|c|c|c|}
    \toprule
$q_J$ & $q_{L}$  & Approach & $J$  & $\tilde\epsilon$\\
    \midrule\midrule
 - & - & Scenario-based  & 2305.55& 0.1248\\
 	\midrule
 1 & 0 & LQ  & 126.44& 1 \\
 0 & 1 & LQ  & 4347.20& 0.9724 \\
 0.2 & 9 & LQ  & 2318.50& 0.9960 \\
    \bottomrule
    \end{tabular}%
  \label{tab:exaPERF}%
\end{table}%

As it can be seen, by means of the scenario-based approach a good trade-off between the cost function $J$ and the violation $\tilde\epsilon$ can be achieved. In particular, though the required value of $\epsilon=0.1$ is not achieved (as it turned out to be infeasible), the actual violation $\tilde\epsilon$ results quite close to the desired one, because $h_1^\star$ and $h_2^\star$ have been properly pushed toward 0 by the the first program \eqref{eq:prob_Sh} in the cascade of problems \eqref{eq:prob_Scenario}, see Fig. \ref{fig:distr}(a).

As for the LQ controller, instead, when only the mass displacements are accounted for ($q_J=1$ $q_{L}=0$) the achieved cost function $J$ is much improved, but, on the other hand, the state constraint is violated by a huge extent as shown in Fig. \ref{fig:distr}(b). When, instead, in the design of the LQ controller only the speeds of the masses are accounted for ($q_J=0$ $q_{L}=1$), the cost function $J$ turns out to be significantly increased with respect to the one obtained by the scenario-based controller. Moreover the speed constraint turns out to be violated with large probability, because the control action tends to excessively reduce the speed at subsequent time steps while it maintains a high speed at the first time instant, see Fig. \ref{fig:distr}(c). In a third design the LQ ($q_J=0.2$ $q_{L}=9$), weights are chosen so as to obtain a performance $J$ similar to the one obtained by the scenario-based controller. Also in this case, however, the probability of constraint violation is high, because the speed constraint is significantly violated in the first time instant, while the masses speeds are excessively reduced in the subsequent time instants, see Fig. \ref{fig:distr}(d).

\begin{figure}[th!]
\centering
\subfloat[Scenario-based controller]{\includegraphics[trim=3.5cm 14cm 4cm 8.5cm, clip=true,scale=0.6]{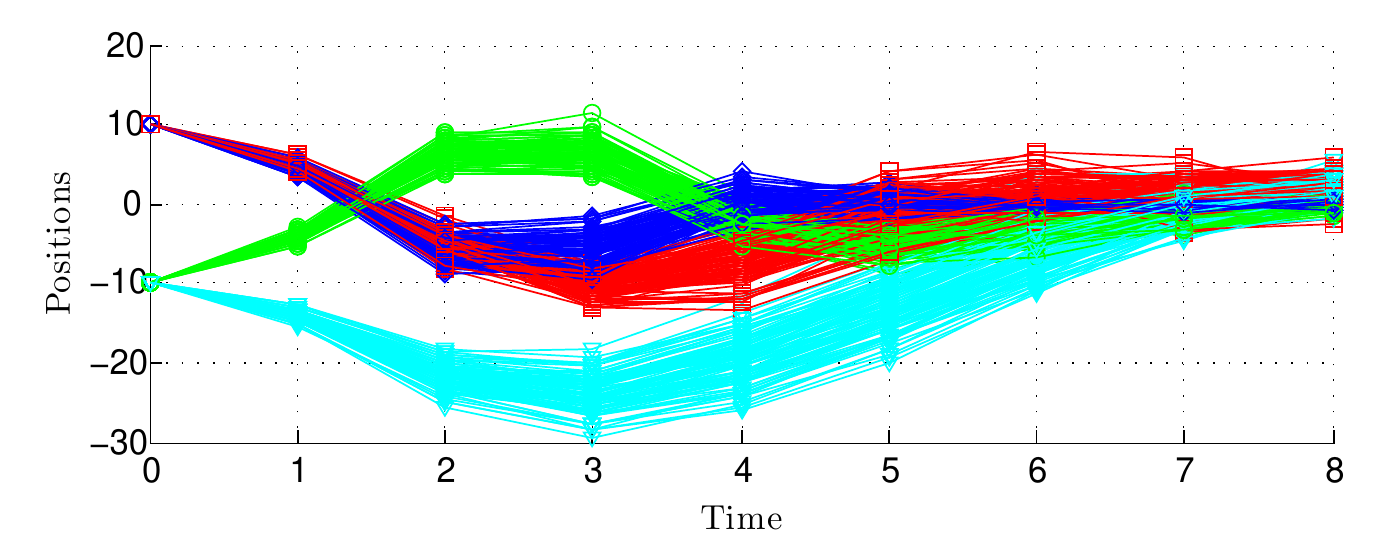}}\\
\subfloat[LQ controller $q_J=0.2$, $q_{L}=9$]{\includegraphics[trim=3.5cm 8.4cm 4cm 14.2cm, clip=true,scale=0.6]{./PaperPosTrJY.pdf}}\\
\caption{Displacements of the masses: $d_1$ (blue diamonds), $d_2$ (green circles), $d_3$ (red squares), $d_4$ (cyan triangles).}
\label{fig:pos}
\end{figure}

\begin{figure}[th!]
\centering
\subfloat[Scenario-based controller]{\includegraphics[trim=3.5cm 14cm 4cm 8.5cm, clip=true,scale=0.6]{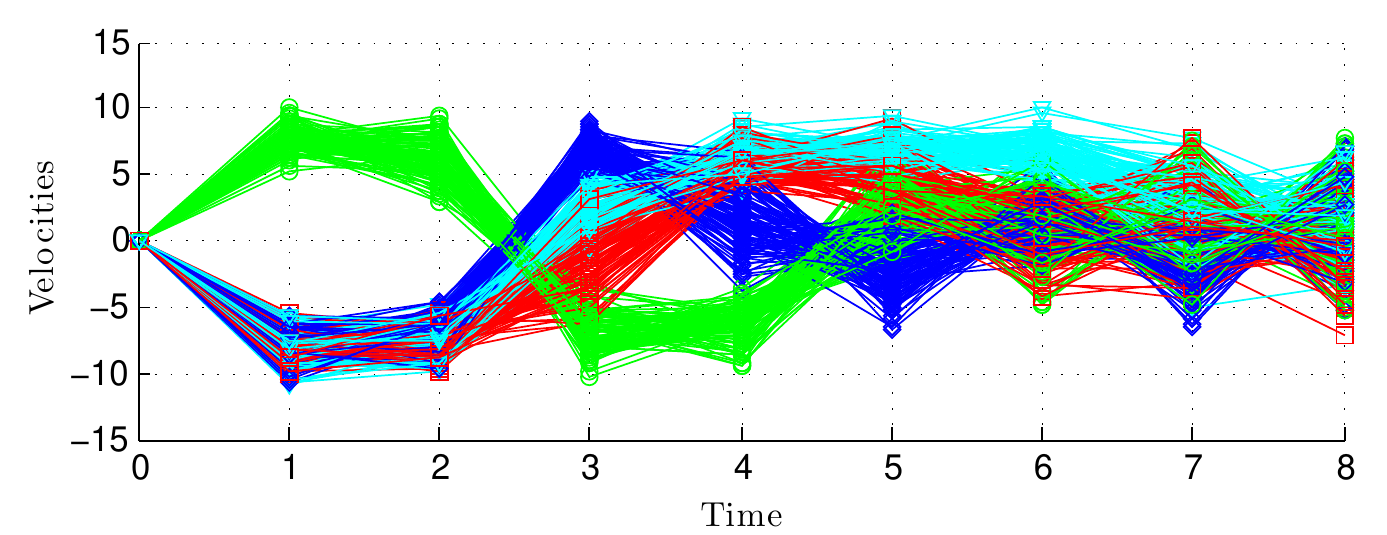}}\\
\subfloat[LQ controller $q_J=0.2$, $q_{L}=9$]{\includegraphics[trim=3.5cm 8.4cm 4cm 14.2cm, clip=true,scale=0.6]{./PaperVelTrJY.pdf}}\\
\caption{Velocities of the masses: $\dot d_1$ (blue diamonds), $\dot d_2$ (green circles), $\dot d_3$ (red squares), $\dot d_4$ (cyan triangles).}
\label{fig:vel}
\end{figure}

The different behaviors of the controllers can be also appreciated by analyzing the state trajectories corresponding to 100 disturbance realizations as depicted in Fig. \ref{fig:pos} (displacements) and Fig. \ref{fig:vel} (velocities). The scenario-based controller exploits the allowed speed to steer the masses close to their nominal position. On the contrary, the LQ control policy leads to the violation of the state constraint in the first time instant, while, in the other instants, the speed is kept conservatively small, and the masses are not steered toward the nominal position.

\bibliographystyle{elsarticle-num}
\bibliography{sMPC_relax_h}

\end{document}